\documentclass[reqno]{amsart}

\usepackage{amsmath, amsfonts, amsthm, amssymb, color,  graphicx, mathrsfs, cite}
\usepackage{comment,stmaryrd}

\theoremstyle{plain}
\newtheorem{theorem}{Theorem}[section]
\newtheorem{lemma}{Lemma}[section]
\newtheorem{proposition}{Proposition}[section]
\newtheorem{corollary}{Corollary}[section]

\theoremstyle{definition}
\newtheorem{definition}{Definition}[section]

\theoremstyle{remark}

\newtheorem{example}{Example}[section]

\numberwithin{equation}{section}
\allowdisplaybreaks

\usepackage{ifpdf}
\ifpdf \usepackage[colorlinks=true, citecolor=blue, linkcolor=blue, urlcolor=blue]{hyperref} \fi

\def\thrm{\begin{theorem}}
\def\thrml#1{\begin{theorem}\label{#1}}
\def\ethrm{\end{theorem}}
\def\lmm{\begin{lemma}}
\def\lmml#1{\begin{lemma}\label{#1}}
\def\elmm{\end{lemma}}
\def\dfntn{\begin{definition}}
\def\dfntnl#1{\begin{definition}\label{#1}}
\def\edfntn{\end{definition}}
\def\crllr{\begin{corollary}}
\def\crllrl#1{\begin{corollary}\label{#1}}
\def\ecrllr{\end{corollary}}
\def\xmpl{\begin{example}}
\def\xmpll#1{\begin{example}\label{#1}}
\def\exmpl{\end{example}}
\def\nmrt{\begin{enumerate}}
\def\enmrt{\end{enumerate}}
\def\qtn{\begin{equation}}
\def\qtnl#1{\begin{equation}\label{#1}}
\def\eqtn{\end{equation}}
\def\prpstn{\begin{proposition}}
\def\prpstnl#1{\begin{proposition}\label{#1}}
\def\eprpstn{\end{proposition}}

\def\tm#1{\item[{\rm (#1)}]}
\def\proof{{\bf Proof}.\ }
\def\eprf{\hfill$\square$}

\DeclareMathOperator{\aut}{Aut}

\DeclareMathOperator{\supp}{supp}

\newcommand{\F}{\mathbb{F}}

\newcommand{\Z}{\mathbb{Z}}

\def\qaq{\quad\text{and}\quad}

\def\lg{\langle}
\def\rg{\rangle}

\def\cZ{\mathcal {Z}}
\def\cT{\mathcal {T}}
\def\cA{\mathcal {A}}

\def\cK{{\mathcal K}}

\def\cD{{\mathcal D}}

\def\qaq{\quad\text{and}\quad}

\begin{document}

\title{Schur rings over Free Abelian Group of Rank Two}
\maketitle

\begin{center}
{\author Gang Chen$^{a}$,  Jiawei He$^{b,*}$ and Zhiman Wu$^{c}$}
\end{center}

\vskip 3mm

\begin{abstract}
Schur rings are a type of subrings of group rings afforded by a partition of the underlined group. In this paper, Schur rings over free abelian group of rank two are classified under the assumption that one of the direct factor is a union of some basic sets. There are eight different types, and all but one type of which are traditional. 
\vskip 2mm

\end{abstract}

\maketitle
\date{}

\maketitle


\renewcommand{\thefootnote}{\empty}
\footnotetext{$^a$ Email: {chengangmath@mail.ccnu.edu.cn}. School of Mathematics and Statistics, and Key Laboratory of Nonlinear Analysis $\&$ Applications (Ministry of Education), Central China Normal University, Wuhan 430079, China}
\footnotetext{ $^{b,*}$ corresponding author,  Email:  hjwywh@mails.ccnu.edu.cn. School of Mathematics and Information Science, Nanchang Hangkong University, Nanchang, China.}
\footnotetext{$^c$  Email: {1793938418@qq.com}. Fengcheng No.9 Middle School, Yichun 331100, China}

\section{Introduction}\label{in}

The concept of Schur ring was first introduced by Schur in \cite{Sch} and later studied by Wielandt in \cite{Wie} in their work concerning permutation groups. Afterwards, of particular importance is the fact that Schur rings have been found to have many applications in algebraic combinatorics. For the applications of Schur rings to association schemes, see for example, subsections 2.4 and 4.4 of \cite{CP}.

\medskip

Classification of Schur rings over a given group  $G$ has been a subject of much research interest.
Most papers have studied Schur rings over finite groups.
Leung and Man have classified Schur rings over finite cyclic groups in \cite{4} and \cite{5}. According to their results, there are four types of Schur rings over finite cyclic groups which are: partitions induced by direct product decomposition of groups ({\it Schur rings of tensor produt type});  partitions induced by orbits of automorphism subgroups ({\it orbit Schur rings}); partitions induced by cosets of normal subgroups ({\it Schur rings of wedge product type}); {\it trivial Schur ring}. We say that a Schur ring over a group is a {\it traditional Schur ring} if it is one of these four types.

\medskip

This paper is a continuation of \cite{BBHMT}.  The authors of the current paper introduced a method for constructing all Schur rings over $\cZ \times \cZ_3$ in \cite{CH} and over infinite dihedral groups in \cite{CH24}. The idea involved is extending a Schur ring over a subgroup to one over the whole group. In this paper, we shall look at another group, that is, the free  abelian group  $\cZ \times \cZ$  of rank two. The structure of  any Schur ring $\cA$ over $\cZ \times \cZ$ which satisfies one of the direct factor $\cZ$ is an $\cA$-subgroup can be completely characterized. There are eight different types, and all but one type of which are traditional.

\medskip

Throughout,  let $G$ denote an arbitrary group with identity element $1$, and $\mathbb{F}$ an arbitrary field of characteristic $0$ and we assume that it contains $\Z$ as a subring. The group algebra over $G$ with coefficients in $\F$ is denoted $\F[G]$.  Let $\Z^{\#}$ denote the set of nonzero integers. For an abelian group $G$, it has an automorphism sending any element to its inverse and the corresponding Schur ring is denoted $\F[G]^{\pm}$. 

\medskip

Let $\cZ$ denote an infinite cyclic group, written multiplicative, and let $\cZ\times \cZ=\langle a\rangle \times\langle b\rangle$ denote a free  abelian group of rank two. The following is the main result of the current paper.

\medskip

\thrml{1451b} Let $\cA$ be a Schur ring over~$\cZ\times\cZ$. If there exists $k\in \Z^{\#}$~such that~$\langle a^k\rangle$~is an ~$\cA$-subgroup.
 Then one the following holds:
\nmrt
\tm{i} $\cA=\F[\cZ\times \cZ]~;$

\tm{ii}  $\cA=\F{[\cZ\times\cZ]}^{\pm}~;$

\tm{iii} $\cA= \F[\cZ]\otimes\F[\cZ]^{\pm}$;

\tm{iv} $\cA=\F[\cZ]^{\pm}\otimes\F[\cZ]^{\pm}$;

\tm{v} $\cA=\mathrm{Span}_{\F}\left\{{\underline D}: D \in \cD(\cA)\} \right\}$, where~$\cD(\cA)=\{\{a^{nj+i}b^{-j},a^ib^{j}\}: i,j\in \Z\}$ for some fixed $n\in\Z^{\#}$;

\tm{vi} $\cA=\mathrm{Span}{_\F}\left\{{\underline D} : D \in \cD(\cA)\} \right\}$, where~$\cD(\cA)=\{\{a^{nj-i}b^{j},a^ib^{j}\}: i,j\in \Z\}$ for some fixed $n\in\Z^{\#}$;

\tm{vii} $\cA=\mathrm{Span}_{\F}\left\{\:{\underline D\: : \:D \in \cD(\cA)\:} \right\}$, where ~$\cD(\cA)=\{\{a^{i}b^{j},a^{nj-i}b^j,a^{-i}b^{-j},a^{i-nj}b^{-j}\} : \:i,j\in \Z\:\}$ for some ~$n\in\Z^{\#}$; 

\tm{viii} $\cA=\mathrm{Span}_{\F}\left\{\:{\underline D\: : \:D \in \cD(\cA)\:} \right\}$, 
where~$\cD(\cA)=\{\{a^i,a^{-i}\}: i\in \Z\}\cup\{\{a^{nj+i}b^{-j},a^ib^{j}\}: \:i,j\in \Z,j\neq0\:\}$ for some fixed $n\in\Z^{\#}$.

\enmrt

Furthermore, all but the last type is traditional.
\ethrm

\section{Preliminary}\label{Pr}
For any finite non-empty subset $D\subseteq G$, denote $\underline{D}=\sum_{g\in C}g$ and call $\underline{D}$ a {\it simple quantity}. In addition, we define $D^*=\{g^{-1}: g\in D\}$. A partition $\cD$ of $G$ is said to be  {\it finite support} provided any $D\in \cD$ is a finite subset of $G$. Our notation on Schur rings is taken from \cite{Po}.

 \subsection{Schur Rings}
\dfntnl{1055a} Let $\cD$ be a finite support partition of $G$ and $\cA$ a subspace of $\mathbb{F}[G]$ spanned by $\{\underline{D}: D\in \cD\}$. We say that $\cA$ is a {\it Schur ring} (or $S$-ring) over $G$ if
 \nmrt
  \tm{1} $\{1\} \in \cD$,
  \tm{2} for any $D\in \cD$, $D^*\in \cD$,
  \tm{3} for all $D_1, D_2\in \cD$, $\underline{D_1}\cdot\underline{D_2}=\sum_{D\in \cD}\lambda_{D_1D_2D}\underline{D}$, where all except finitely many $\lambda_{D_1D_2D}$ are equal to $0$.
 \enmrt
\edfntn

For a Schur ring $\cA$ over $G$, the associated partition $\cD$ is denoted $\cD(\cA)$ and each element in $\cD(\cA)$ is called a {\it basic  set}
of $\cA$.  Each $\lambda_{D_1D_2D}$ is called a {\it structure constant}. We say that a subset $D$ of $G$ is an {\it $\cA$-set} if $D$ is a union of some basic sets of $\cA$. An {\it $\cA$-subgroup} is simultaneously an $\cA$-set and a subgroup.

\medskip

Suppose $\alpha=\sum_{g\in G}\alpha_gg\in \F[G]$ and $\beta=\sum_{g\in G}\beta_gg\in \F[G]$, define
$$
\alpha^*=\sum_{g\in G}\alpha_gg^{-1}~ \mathrm{and }~\supp(\alpha)=\{g\in G: \alpha_g\ne 0 \}.
$$
 Additionally,  the {\it Hadamard product} of $\alpha$ and $\beta$ is defined as 
 $$
 \alpha\circ \beta=\sum_{g\in G}\alpha_g\beta_gg.
 $$

\medskip

The following lemma is well-known, which is \cite[Lemma 2.15]{BBHMT}.

\medskip

\lmml{1718c8} Let $\cA$ be a Schur ring over a group $G$. Suppose $C=\{g\}$ is a basic set. Then for any basic set $D$, $gD$  and $Dg$ are basic sets.

\elmm

The following theorem was proved by Wielandt in \cite{Wie2} when $G$ is finite. The general case appeared as \cite[Corollary 2.12]{BBHMT}.

\thrml{1505b} Let $G$ be a group and $\cA$ a subring of $\F[G]$ which is spanned by $\{\underline{D}: D\in \cD(\cA)\}$ as an $\F$-vector space. Then $\cA$ is a Schur ring over $G$ with basic sets from $\cD(\cA)$ if and only if $\cA$ is closed under $\circ$ and $*$, $1\in \cD(\cA)$, and for all $g\in G$ there exists some $D\in \cD(\cA)$ such that $g\in D$.
\ethrm

Let $f: \F\rightarrow \F$ be a function such that $f(0)=0$. For any $\alpha=\sum_{g\in G}\alpha_gg$, we set
$$
f(\alpha)=\sum_{g\in G}f(\alpha_g)g.
$$

The following proposition was proved by Wielandt in \cite[Proposition 22.3]{Wie2}. The general case appeared as \cite[Propostion 2.4]{BBHMT}.

\prpstnl{913a}Let $\cA$ be a Schur ring over a group $G$ and the function be as above. Then
$f(\alpha)\in \cA$ whenever $\alpha\in \cA$.
\eprpstn

Let $f$ take value $1$ at one non-zero number and $0$ at other numbers.  We get the following statement, which is known as the Schur-Wielandt principle; see \cite[Corollary 1.10]{Po}.

\crllrl{10131c} Let $\cA$ be a Schur ring over $G$. For any $\alpha=\sum_{g\in G}\alpha_gg\in \cA$, suppose $g\in \supp(\alpha)$ with $\alpha_g=c\ne 0$. Then the set
$$
\{g\in G: \alpha_g=c\}
$$
is an $\cA$-set.
\ecrllr

The following proposition, which was first proved by Wielandt in \cite{Wie2} and was generalized in \cite{BBHMT},  tells us how to generate $\cA$-subgroups in a Schur ring $\cA$.

\prpstnl{956b} Let $\cA$ be a  Schur ring over $G$. Let $\alpha\in \cA$ and $H=\lg \supp(\alpha)\rg$. Then~$H$ is an $\cA$-subgroup.

\eprpstn

Suppose that $H$ is an $\cA$-subgroup, then $\{D: D\in \cD(\cA), D\subseteq H\}$
consists of  a basic set of a Schur ring over $H$. It is denoted $\cA_H$.

\medskip

When we have normal $\cA$-subgroup, we can construct  Schur ring over the factor group as shown in the following lemma.

 \lmml{1440a}(\cite[Lemma 1.2]{LM}) Let $\varphi: G\rightarrow H$ be a group homomorphism with kernel $K$ and $\cA$  a Schur ring over $G$. Suppose that $K$ is an $\cA$-subgroup. Then the image $\varphi(\cA)$ is a Schur ring over $\varphi(G)$, where $\cD(\varphi(\cA))=\{\varphi(C): C\in \cD(\cA)\}$.
\elmm

In particular, if $H=G/K$, the factor group of $G$ over $K$, and $\varphi$ is the natural homomorphism,  the corresponding Schur ring is denoted $\cA_{G/K}$.

\medskip

\subsection{Traditional Schur rings}
\
\medskip

\textbf{Discrete and trivial Schur rings}.
The group ring $\F[G]$ itself forms a Schur ring over $G$, which is called
the {\it discrete Schur ring} over $G$. When $G$ is a finite group, the trivial partition $\{1,G\setminus 1 \}$ produces
a Schur ring which is known as the {\it trivial Schur ring} over $G$.

\medskip
\textbf{Tensor products}. Assume that $G =G_1 \times G_2 $, $\cA_1$ and $\cA_2$ are Schur rings over ${G_1}$ and ${G_2}$, respectively. Then the
set 
$$
\{CD : C\in \cD(\cA_1), D\in \cD(\cA_2)\}
$$
consists of the basic sets of a Schur ring over $G$, which was first proved by Wielandt in \cite{Wie} when group $G$ is a finite group. If $G$ is an infinite group, it is easy to see that this conclusion also
holds.  The corresponding Schur ring is called the {\it tensor product} of $\cA_1$ and $\cA_2$, denoted $\cA_1\otimes \cA_2$. 

\medskip

\textbf{Orbit Schur rings}. If $\cK$ is a finite subgroup of $\aut(G)$, the set of elements of $\F[G]$ fixed by $\cK$ is a Schur ring over $G$, denoted $\F[G]^{\cK}$ and called the {\it orbit Schur ring} associated with $\cK$. Obviously, the group ring $\F[G]$ itself also forms an orbit Schur ring {associated with the trivial automorphism group}.

\medskip

\textbf{Wedge products}. A Schur ring $\cA$ over $G$ is a {\it wedge product} if there exist nontrivial proper $\cA$-subgroups $H, K$ such that $K \leq H, K \unlhd G$, and every basic set outside $H$ is a union of $K$-cosets. In this case, the series
$$
1<K \leq H<G
$$
is called a wedge-decomposition of $\cA$. Furthermore, we write $\cA=\cA_1 \wedge \cA_2$, where $\cA_1=\cA_H$ and $\cA_2=\cA_{G / K}$.

\medskip

A Schur ring over a group $G$ is called {\it traditional} if it is either a trivial Schur ring (when $G$ is finite), or an orbit Schur ring, or a tensor product of Schur rings over proper subgroups, or a wedge product. A group $G$ is said to be {\it traditional} if each Schur ring over it is traditional.

\medskip

\thrml{1514a}(\cite[Theorem  3.3]{BBHMT}) The only Schur rings over the infinite cyclic group $\cZ:=\lg z \rg$ are either discrete or symmetric Schur ring. Here each basic set of the discrete Schur ring is a singleton set and each basic set of the  the symmetric Schur rings has the form $\{z^i, z^{-i}\}$ for some nonnegative integer $i$. Note that the latter Schur ring is $\F[\cZ]^{\pm}$. 
\ethrm

\medskip

\section{Proof of the main Results}

For $n \in \Z$, define the $n$th Frobenius map of $\F[G]$ as
$$
\alpha=\sum_{g \in G} \alpha_g g \mapsto \alpha^{(n)}=\sum_{g \in G} \alpha_g g^n,
$$
for any $\alpha=\sum_{g\in G}\alpha_gg\in \F[G]$. 
\medskip

\lmml{1718c111} (\cite[Theorem 2.20]{BBHMT}) Suppose that $\cA$ is a Schur ring over $G$ where $G$ is an abelian group. Let $m$ be an integer coprime to the orders of all torsion elements of $G$. Then, for all $\alpha\in \cA$ we have $\alpha^{(m)} \in \cA$. In particular, if $G$ is torsion-free, then, all Schur rings over $G$ are closed under all Frobenius maps.
\elmm
\medskip

\lmml{1718c112} (\cite[Lemma 2.21]{BBHMT}) Let $G$ be an abelian group such that the torsion subgroup $T(G)$ has finite exponent. Let $\cA$ be a Schur ring over $G$. Then, $T(G)$ is an $\cA$-subgroup.
\elmm

\medskip
{\bf Proof of Theorem \ref{1451b}}\,
In the sequel, we fix the following notation:
$$
G=\cZ\times\cZ=\lg a\rg\times \lg b\rg.
$$ 
Let
$$
\varphi:G\rightarrow G/\lg a^k\rg
$$
be the natural
homomorphism, and set~
$$
\cT:=\varphi(\cA).
$$
Then~$\cT$~is~a Schur ring over $G/\lg a^k\rg$~by Lemma \ref{1440a}. Moreover, according to Lemma \ref{1718c112},
 $$
 T(\varphi(G))=\lg \varphi(a)\rg
 $$
  is a ~$\cT$-subgroup.  It follows that
  $$
  \lg a\rg=\varphi^{-1}(\lg \varphi(a)\rg)
  $$
  is  an~$\cA$-subgroup.
Therefore we may assume that~$\cZ=\lg a\rg$~is an $\cA$-subgroup in the following.

\medskip

Choose a basic set $C$ such that  $a\in C$. Then~
$$
C=\{a\}~\mathrm{or} ~C=\{a,a^{-1}\}
$$
by Theorem \ref{1514a}. Hence $\{a^i,a^{-i}\}$ are~$\cA$-sets for all~$i\in \Z$. In addition, as $G/\lg a\rg$~is cyclic,  the basic set containing $\varphi(b)$ is equal to
$$
\{\varphi(b)\}~ \mathrm{or}~\{\varphi(b),\varphi(b^{-1})\}
$$
 by Theorem \ref{1514a}.
 \medskip

\thrml{1718c} All conditions are the same as in Theorem 1.1, and let $D\in \cD(\cA)$ with $b\in D$. Then one of the following holds:
 \nmrt

\tm{i} $D=\{b\}$;

 \tm{ii} $D=\{b,b^{-1}\}$;

 \tm{iii} $D=\{ba^{i_0},b\}$ for some ~$i_0\in\Z^{\#}$;

\tm{iv} $D=\{b^{-1}a^{i_1},b\}$ for some ~$i_1\in\Z^{\#}$;

\tm{v} $D=\{b, ba^{i_2},b^{-1}, b^{-1}a^{-i_2}\}$ for some ~$i_2\in \Z^{\#}$.

\enmrt
\ethrm

\proof Observe that
$$
D\subseteq\varphi^{-1}(\{\varphi(b),\varphi(b^{-1})\})=b\lg a\rg\cup b^{-1}\lg a\rg.
$$
\medskip

To prove this theorem, we next consider the following two cases.

 \medskip

{\bf Case 1.}\, Suppose
$$
D\subseteq b\lg a\rg.
$$
Set $D:=b\{a^{i_1}, \ldots, a^{i_k}\}$, where $i_1<\cdots<i_k$ and one of the $i_j$ is zero.

\medskip

Then either $D=\{b\}$ or $k$ is greater than or equal to $2$. If $k>2$, then
\qtnl{2027c}
\{a^{i_1}b,a^{i_k}b\}=D\cap\{a^{i_k-i_1},a^{i_1-i_k}\}D
 \eqtn
 is an $\cA$-set, a contradiction. We conclude that $k=2$ and $D=\{a^{i_0}b,b\}$ for some $i_0\in\Z^{\#}$ in this case.

\medskip

{\bf Case 2.}\, Now suppose
$$
D=b\{a^{i_1}, \ldots, a^{i_k}\}\cup b^{-1}\{a^{j_1}, \ldots, a^{j_t}\},
$$
where $i_1<\cdots <i_k$ and $j_1<\cdots <j_t$.

\medskip

{\bf Claim 1.} \,  {\it We have $k\le 2$ and $t\le 2$}.

\medskip

\proof
If the claim is false, assume on the contrary that $k>2$. Then similar to equality \eqref{2027c}, one can see that $D\cap\{a^{i_k-i_1},a^{i_1-i_k}\}D$   contains $\{ba^{i_1}, ba^{i_k}\}$ but not $ba^{i_3}$, which is a contradiction.  Thus, $k\le 2$ and similar argument shows that $t\le 2$. The claim follows.

\eprf

{\bf Claim 2.}\, {\it We have $k=t$.}

\medskip

\proof If the claim is false,  we first assume that $k=1, t=2$. Let
$$
D=\{b\}\cup b^{-1}\{a^{j_1}, a^{j_2}\},
$$
where $j_1<j_2$.
\medskip

Then  $\cA$ contains the following element:
\qtnl{2122a}
c:=\underline{D}(a^{-j_1}+a^{j_1})=b(a^{-j_1}+a^{j_1})+b^{-1}(1+a^{j_2-j_1}+a^{2j_1}+a^{j_1+j_2}).
\eqtn
Observe that
$$
D^*=\{b^{-1}, ba^{-j_1}, ba^{-j_2}\}.
$$
Thus, $D^*$ is contained in $\supp(c)$ and this fact forces $j_1=-j_2$. It then follows that 
$$
c:=b(a^{-j_1}+a^{j_1})+2b^{-1}+b^{-1}(a^{-2j_1}+a^{2j_1}).
$$
Thus, $c=2\underline{D^*}$, which is absurd.

\medskip

Now assume that $k=2$ and $t=1$.  Let
$$
D=\{b, ba^{i_1},  b^{-1}a^{j_1}\}.
$$
Then
$$
c':=\underline{D}(a^{-j_1}+a^{j_1})=ba^{-j_1}+ba^{j_1}+ba^{i_1-j_1}+ba^{i_1+j_1}+b^{-1}+b^{-1}a^{2j_1}\in \cA.
$$
Note that 
$$
D^*=\{b^{-1}, b^{-1}a^{-i_1}, ba^{-j_1}\}. 
$$
One can see that $D^*\subseteq \supp(c')$. This implies that $i_1=-2j_1$. It then follows that 
$$
c'=2ba^{-j_1}+b^{-1}+b^{-1}a^{2j_1}+ba^{j_1}+ba^{-3j_1}. 
$$ 
This yields that $c'=2\underline{D^*}$, which is impossible. 
	
\eprf

\medskip

Now assume that $k=t=2$. Let
$$
D=\{b, ba^{i_2}, b^{-1}a^{j_1}, b^{-1}a^{j_2}\},
$$
where $i_2\in \Z^{\#}$ and $j_1<j_2$.
Then $\cA$ contains the following element:
\begin{displaymath}
	\begin{split}
&\underline{D}(a^{i_2}+a^{-i_2})=b+ba^{i_2}+ba^{-i_2}+ba^{2i_2}\\
+&b^{-1}a^{j_1+i_2}+b^{-1}a^{j_1-i_2}+b^{-1}b^{j_2+i_2}+b^{-1}a^{j_2-i_2}.
\end{split}
\end{displaymath}
One can see that  $j_2+i_2=j_1$ or $j_1+i_2=j_2$.

\medskip

Consequently, there exist  integers $i_2\ne 0$ and $j_1$ such that
$$
D=\{b, ba^{i_2},b^{-1}a^{j_1},b^{-1}a^{i_2+j_1}\}.
$$
Then
$$
D^{*}=\{b^{-1}, b^{-1}a^{-i_2}, ba^{-j_1}, ba^{-i_2-j_1}\}
$$
is also a basic set. Note that
  \begin{displaymath}
  \begin{split}
  d:=&(a^{j_1}+a^{-j_1})\underline{D}\\
  =&ba^{j_1}+ba^{-j_1}+ba^{i_2+j_1}+ba^{i_2-j_1}+\\
  &b^{-1}a^{2j_1}+b^{-1}+b^{-1}a^{i_2+2j_1}+b^{-1}a^{i_2}\in \cA.
  \end{split}
  \end{displaymath}
Thus,  $D^{*}\subseteq \supp(d)$. Hence
$$
2j_1=-i_2, \,\, {\rm or}\, \, j_1=-i_2.
$$

\medskip

{\bf Claim 3.}\, {\it The case where $2j_1=-i_2$ can not happen.}
\medskip

\proof Towards a contradiction, assume that $2j_1=-i_2$. Then
$$
D=\{b, ba^{-2j_1}, b^{-1}a^{j_1}, b^{-1}a^{-j_1}\}, \qaq D^*=\{b^{-1}, b^{-1}a^{2j_1}, ba^{-j_1}, ba^{j_1}\}.
$$
In addition,
\begin{displaymath}
  \begin{split}
  f:=&(a^{j_1}+a^{-j_1})\underline{D}\\
  =&2(b^{-1}+ba^{-j_1})+ba^{j_1}+ba^{-3j_1}+b^{-1}a^{2j_1}+b^{-1}a^{-2j_1}\in \cA.
  \end{split}
  \end{displaymath}
Thus, $f=2\underline{D^*}$, which forces $j_1=0$ and thus $i_2=0$, a contradiction.

\medskip

 Therefore,  $j_1=-i_2$. It then follows that
 $$
 D=\{b,ba^{i_2}, b^{-1}, b^{-1}a^{-i_2}\}
 $$
 for some $i_2\in \Z^{\#}$.

 \medskip
 Finally, if $k=t=1$, then statement $(iv)$ in the theorem happens. The proof of this theorem is complete.

\eprf

\crllrl{959b} All conditions are the same as in Theorem 1.1, and let $D\in \cD(\cA)$ with $b\in D$. If there exists $i_0,j_0\in \Z,j_0\neq0$ such that ~$\{a^{i_0}b^{j_0}\}\in \cD(\cA)$, then one of the following holds:
 \medskip

 \nmrt

\tm{i} $D=\{b\}$;

\tm{ii} $D=\{a^{n_0}b,b\}$, where the integer $n_0$  satisfies  $2i_0=n_0j_0$;

\tm{iii} $D=\{a^{n_1}b,a^{-n_1}b^{-1},b,b^{-1}\}$,  where the integer $n_1$  satisfies $2i_0=n_1j_0$.

\enmrt

\ecrllr

\proof
By Theorem \ref{1718c8},
$$
\{a^{2i_0}b^{j_0},b^{j_0}\}=\{a^{i_0},a^{-i_0}\}\{a^{i_0}b^{j_0}\}
$$
 is an $\cA$-set. Also, $D^{(j_0)}$~is an
~$\cA$-set by  Lemma \ref{1718c111}.
Thus, either
$$
\{b^{j_0}\}=D^{(j_0)}\bigcap\{a^{2i_0}b^{j_0},b^{j_0}\}
$$
is an $\cA$-set, or
$$
\{a^{2i_0}b^{j_0},b^{j_0}\}\subseteq D^{(j_0)},  i_0\neq 0.
$$

If~$\{b^{j_0}\}$~is a basic set. Then $\lg b\rg$~is an $\cA$-subgroup as what we have seen at the beginning of the proof. This implies that 
$$
D=\{b\}~\mathrm{or}~D=\{b,b^{-1}\}.
$$
By way of contradiction, assume that~$D=\{b,b^{-1}\}$, then $D^{(j_0)}=\{b^{j_0},b^{-j_0}\}$~ is a basic set by Theorem \ref{1514a}, a contradiction.  Hence, in this case,~$D=\{b\}$. 

\medskip

Thus, we may assume that~$\{a^{2i_0}b^{j_0},b^{j_0}\}\subseteq D^{(j_0)}$ with $i_0\neq 0$.~Then, by Theorem \ref{1718c}, we have the following two cases:

\medskip

The first case : $D=\{a^{n_0}b,b\}$,  where the integer $n_0$  satisfies  $2i_0=n_0j_0$;

\medskip

The second case: ~$D=\{a^{n_1}b,a^{-n_1}b^{-1},b,b^{-1}\}$, where the integer $n_1$  satisfies $2i_0=n_1j_0$.

\medskip

The proof of this corollary is complete.

\eprf
\medskip

Now, we return to the proof of the main theorem. By Theorem \ref{1514a}, we have the following two cases:

\medskip

{\bf Case 1.}\, Suppose that $C=\{a\}$. Then ~$\{a^i\}\in \cD(\cA)$ for ~$i\in \Z$. For ~$n\in\Z^{\#}$, observe that
$$
\{a^{-n}\}\{a^{n}b,b\}=\{b,a^{-n}b\}\neq \{a^{n}b,b\},
$$
and that
\begin{displaymath}
\begin{split}
&\{a^{-n}\}\{a^{n}b,a^{-n}b^{-1},b,b^{-1}\}\\
=&\{b,a^{-2n}b^{-1},a^{-n}b,a^{-n}b^{-1}\}\\
\neq &\{a^{n}b,a^{-n}b^{-1},b,b^{-1}\}.
\end{split}
\end{displaymath}
 By Theorem \ref{1718c}, we conclude that $$
 D=\{b\}, ~\mathrm{or}~D=\{b,b^{-1}\}, ~\mathrm{or}~D=\{a^{n_1}b^{-1},b\},
 $$
 for some  $n_1\in\Z^{\#}$.
  \medskip

 {\bf Subcase 1.1.}\, Suppose $D=\{b\}$. Then ~$\{b^j\}$~is a basic set for all $j\in \Z$.~ So, by Lemma \ref{1718c8},~ ~$\{a^ib^j\}$~are basic sets for all ~$i,j\in \Z$. Hence,
 $$
 \cA=\F[\cZ\times\cZ].
 $$
 This is case $(i)$ in Theorem \ref{1451b}.

 \medskip

 {\bf Subcase 1.2.}\, Assume $D=\{b,b^{-1}\}$. Then by Lemma \ref{1718c111} and Theorem 2.2, ~$\{b^j,b^{-j}\}$~is a basic set for all~$j\in \Z$. Furthermore,  ~$\{a^ib^j,a^ib^{-j}\}$~are  basic sets for ~$i,j\in \Z$. Thus,
$$
\cA=\F[\cZ]\otimes \F[\cZ]^{\pm},
$$
which is  case $(iii)$ in Theorem \ref{1451b}.

\medskip

 {\bf Subcase 1.3.}\, Assume that~$D=\{a^nb^{-1},b\}$, ~where~$n\in \Z^{\#}$. According to Lemma \ref{1718c111}, it follows that~$\{a^{nj}b^{-j},b^{j}\}$~is a basic set for each~$j\in \mathbb{Z}$. Therefore
 $$
 \{a^{nj+i}b^{-j},a^ib^{j}\}
 $$
 is a basic set for each $i,j\in \Z$.
 \medskip

 Observe that $G$ has an automorphism $\varphi_1$ such that
 $$
 \varphi_1(a)=a,~\varphi_1(b)=a^{n}b^{-1}.
 $$
 Also, it is obvious that the $\lg \varphi_1\rg$-orbit contains $a^ib^j$ is $\{a^{nj+i}b^{-j},a^ib^{j}\}$.

 \medskip

This means that
$$
\cA=\F[G]^{\lg \varphi_{1}\rg},
$$
i.e.,$\cA$ is the orbit Schur ring associated with $\lg \varphi_1\rg$.  This is case $(v)$ in Theorem \ref{1451b}.

\medskip

 {\bf Case 2.}\, Suppose~$C=\{a,a^{-1}\}$.  Then ~$\{a^i,a^{-i}\}$~is a basic set for any ~$i\in\Z$.

  \medskip

 {\bf Subcase 2.1.}\, Suppose that $D=\{b\}$. Hence, for all $m\in \mathbb{Z}\setminus \{0\}, n\in \mathbb{Z}$,
 $$
 \{a^m, a^{-m}\}\in \cD(\cA), ~\mathrm{and}~\{a^mb^n, a^{-m}b^n\}\in \cD(\cA).
 $$
 Thus,
$$
\cA=\F[\cZ]\otimes {\F[\cZ]}^{\pm},
$$
which is  case $(iii)$ in Theorem \ref{1451b}.
\medskip

 {\bf Subcase 2.2.}\,Suppose~$D=\{b,b^{-1}\}$. For any nonzero integers $i,j$, $\{b^j,b^{-j}\}$ is a basic set, and hence
 $$
 \{a^ib^j,a^ib^{-j},a^{-i}b^j,a^{-i}b^{-j}\}=\{a^i,a^{-i}\}\{b^j,b^{-j}\}
 $$
 is an $\cA$-set. By Corollary \ref{959b}, we know that $\{a^ib^j\}$~is not a basic set.

 \medskip

If~$\{a^ib^j,a^ib^{-j}\}$~is a basic set for nonzero integers $i$ and $j$, ~then
$$
(a^ib^j+a^ib^{-j})^2=a^{2i}b^{2j}+a^{2i}b^{-2j}+2a^{2i}\in \cA.
$$
By Proposition \ref{913a}, ~$\{a^{2i}\}$~is a basic set, which contradicts the fact that $\{a^{2i},a^{-2i}\}$~is a basic set. It follows that~$\{a^ib^j,a^ib^{-j}\}$ is not a basic set. Similarly, $\{a^ib^j,a^{-i}b^{j}\}$~is also not a basic set.

\medskip

Applying Corollary \ref{959b} again, one can see that the basic set containing ~$a^ib^j$ is either
$$
\{a^ib^j,a^{-i}b^{-j}\}
$$
or
$$
\{a^ib^j,a^ib^{-j},a^{-i}b^j,a^{-i}b^{-j}\}.
$$

Suppose there exist nonzero integers $i_1,i_2,j_1,j_2$ such that
\qtnl{1047a}
\{a^{i_1}b^{j_1},a^{-i_1}b^{-j_1}\},~\{a^{i_2}b^{j_2},a^{i_2}b^{-{j_2}},a^{-{i_2}}b^{j_2},a^{-{i_2}}b^{-{j_2}}\}\in \cD(\cA).
\eqtn
In addition,
{\begin{displaymath}
	\begin{split}
g:=&(a^{i_1-i_2}+a^{i_2-i_1})(b^{j_1-j_2}+b^{j_2-j_1})\\
&=a^{i_1-i_2}b^{j_1-j_2}+a^{i_2-i_1}b^{j_1-j_2}+a^{i_1-i_2}b^{j_2-j_1}+a^{i_2-i_1}b^{j_2-j_1}
\end{split}
\end{displaymath}}
belongs to $\cA$. It follows that
\begin{displaymath}
\begin{split}
h:=&g(a^{i_1}b^{j_1}+a^{-i_1}b^{-j_1})\\
=&a^{2i_1-i_2}b^{2j_1-j_2}+a^{i_2}b^{2j_1-j_2}+a^{2i_1-i_2}b^{j_2}+a^{i_2}b^{j_2}+\\
&a^{-i_2}b^{-j_2}+a^{i_2-2i_1}b^{-j_2}+a^{-i_2}b^{j_2-2j_1}+a^{i_2-2i_1}b^{j_2-2j_1}\in \cA.
\end{split}
\end{displaymath}

One can easily see that the second basic set in \eqref{1047a} is contained in $\supp(h)$. One can easily deduce that at least one of the parameters $i_1, i_2, j_1, j_2$ must be $0$, which is a contradicion.

\medskip

Thus,
$$
\{a^{i}b^{j},a^{-i}b^{-j}\}\in \cD(\cA)
$$ for all ~$i, j\in \Z^{\#}$
or
$$
\{a^{i_2}b^{j_2},a^{i_2}b^{-{j_2}},a^{-{i_2}}b^{j_2},a^{-{i_2}}b^{-{j_2}}\}\in \cD(\cA)
$$
for all $i, j\in \Z^{\#}$.

\medskip

In the first case, one can easily see that
$$
\cA=\F{[\cZ\times\cZ]}^{\pm},
$$
which is case $(ii)$ in Theorem \ref{1451b}.

\medskip

In the second case, one can see that
$$
\cA={\F[\cZ]}^{\pm}\otimes{\F[\cZ]}^{\pm},
$$
which is  case $(iv)$ in Theorem \ref{1451b}.

\medskip

 {\bf Subcase 2.3.}\, Suppose~$D=\{a^nb,b\}$, for some $n\in \Z^{\#}$.~Then, for~$i,j\in\Z,~j\neq0$,~$\{a^{nj}b^j,b^{j}\}$ is an~$\cA$-set. Moreover,
$$
k:=(a^i+a^{-i})(a^{nj}b^j+b^{j})=a^{i+nj}b^j+a^ib^{j}+a^{nj-i}b^j+a^{-i}b^{j}\in\cA
$$

\medskip

Suppose $2i=nj$, then
$$
k=2a^{i}b^{j}+a^{3i}b^{j}+a^{-i}b^{j}\in \cA.
$$
This implies that $\{a^ib^j\}$ is a basic set by Corollary \ref{10131c}.

\medskip

Suppose $2i\neq nj$. As
$$
l:=(a^{i-nj}+a^{nj-i})(a^{nj}b^j+b^{j})=a^{i}b^j+a^{i-nj}b^{j}+a^{2nj-i}b^j+a^{nj-i}b^{j}\in \cA,
$$
it then follows that 
$$
a^ib^{j}+a^{nj-i}b^j=k\circ l\in \cA. 
$$
We then conclude that $\{a^ib^{j},a^{nj-i}b^j\}$ is an $\cA$-set.

\medskip

 Moreover, $\{a^ib^{j},a^{nj-i}b^j\}$  must be a basic set. Otherwise, both $\{a^ib^{j}\}$ and $\{a^{nj-i}b^j\}$ are basic sets. This implies that $\{a^{2i-nj}\}$ is a basic set with $2i-nj\ne 0$, which contradicts the fact that $\{a^{2i-nj}, a^{nj-2i}\}$ is a basic set.

\medskip

As a consequence,~$\{a^ib^{j},a^{nj-i}b^j\}$~is a basic set for any~$i,j\in\Z$. Observe that $G$ has an automorphism $\varphi_2$ such that
$$
\varphi_{2}(a)=a^{-1},~\varphi_{2}(b)=a^{n}b.
$$
Observe that the $\lg \varphi_2\rg$-orbit containing $a^ib^j$ is $\{a^ib^{j},a^{nj-i}b^j\}$.
\medskip

Thus,
$$
\cA=\F[G]^{\lg \varphi_2\rg},
$$
i.e., $\cA$ is the orbit Schur ring associated with $\lg \varphi_2\rg$. This  case $(vi)$  in Theorem \ref{1451b}.

\medskip

 {\bf Subcase 2.4.}\, Suppose that~$D=\{a^nb^{-1},b\}$ for some~$n\in \Z^{\#}$.~Thus,~$\{a^{nj}b^{-j},b^{j}\}$ is an~$\cA$-set for any~$i,j\in\Z,j\neq0$.

\medskip

Moreover,
$$
a^{i+nj}b^{-j}+a^ib^{j}+a^{nj-i}b^{-j}+a^{-i}b^{j}=(a^i+a^{-i})(a^{nj}b^{-j}+b^{j}),
$$
and
$$
a^ib^j+a^{i+nj}b^{-j}+a^{-i-2nj}b^{j}+a^{-i-nj}b^{-j}=(a^{i+nj}+a^{-i-nj})(a^{-nj}b^{j}+b^{-j})
$$
are two elements in $\cA$. So their common support
$\{a^ib^j,a^{i+nj}b^{-j}\}$ is an $\cA$-set. Moreover, it must be a basic set by Corollary 3.1.

\medskip

For $i_1,i_2,j_1,j_2\in \Z,j_1\ne 0 \ne j_2$, we have
\begin{displaymath}
\begin{split}
&(a^{i_1}b^{j_1}+a^{i_1+nj_1}b^{-j_1})(a^{i_2}b^{j_2}+a^{i_2+nj_2}b^{-j_2})\\
=&(a^{i_1+i_2}b^{j_1+j_2}+a^{(i_1+i_2)+n(j_1+j_2)}b^{-j_1-j_2})+
(a^{i_1+i_2+nj_2}b^{j_1-j_2}+a^{i_1+i_2+nj_1}b^{j_2-j_1})\in \cA,
\end{split}
\end{displaymath}
since $i_1+i_2+nj_1=(i_1+i_2+nj_2)+n(j_1-j_2)$.

\medskip

In addition,
\begin{displaymath}
\begin{split}
&(a^{i_1}+a^{-i_1})(a^{i_2}b^{j_2}+a^{i_2+nj_2}b^{-j_2})\\
=&(a^{i_1+i_2+nj_2}b^{-j_2}+a^{i_2+i_1}b^{j_2})+
(a^{i_2-i_1+nj_2}b^{-j_2}+a^{i_2-i_1}b^{j_2})\in \cA.
\end{split}
\end{displaymath}
Thus,
 $$
 \cD(\cA)=\{~\{a^i,a^{-i}\}~|\:i\in \Z\:\}\cup\{\{a^ib^{j},a^{nj+i}b^{-j}\}|\:i,j\in \Z,j\neq0\:\}
 $$
 for some $n\in \Z^{\#}$. This is case $(viii)$ in Theorem \ref{1451b}.

 \medskip

 {\bf Subcase 2.5.}\,Suppose~$D=\{a^{n}b,a^{-n}b^{-1},b,b^{-1}\}$
  for some~$n\in \Z^{\#}$. 
  
  \medskip
  
 Set $u:=\underline{D^{(j)}}$ for nonzero integer $j$. Then,
\begin{displaymath}
\begin{split}
&x:=u(a^i+a^{-i})\\
=&a^{i+nj}b^{j}+a^{i-nj}b^{-j}+a^{i}b^{j}+a^{i}b^{-j}+\\
&a^{nj-i}b^j+a^{-i-nj}b^{-j}+a^{-i}b^{j}+a^{-i}b^{-j}\in\cA,
\end{split}
\end{displaymath}
for all~$i,j\in\Z$ with $j\ne 0$.

\medskip

If  $2i=nj$, we have
$$
x=a^{3i}b^{j}+a^{-3i}b^{-j}+a^{i}b^{-j}+a^{-i}b^{j}+2a^{i}b^{j}+2a^{-i}b^{-j}. 
$$
Then, according to  Proposition \ref{913a},
$$
\{a^{i}b^{j},a^{-i}b^{-j}\}\in \cD(\cA),
$$
for any nonzero integers $i,j$ such that $2i=nj$.

\medskip

If $2i\ne nj$, then
\begin{displaymath}
\begin{split}
&y:=u(a^{i-nj}+a^{nj-i})\\
=&a^{i}b^{j}+a^{i-nj}b^{j}+a^{i-2nj}b^{-j}+a^{i-nj}b^{-j}+\\
&a^{2nj-i}b^{j}+ a^{nj-i}b^{j}+a^{-i}b^{-j}+a^{nj-i}b^{-j}\in\mathcal{A}.
\end{split}
\end{displaymath}
Thus, $$
x\circ y=a^{i}b^{j}+a^{nj-i}b^j+a^{-i}b^{-j}+a^{i-nj}b^{-j}\in \cA.
$$

Suppose $\{a^ib^j\}$ is a basic set. Then
$$
z:=a^ib^j(a^i+a^{-i})=a^{2i}b^j+b^j\in \cA.
$$
Thus,
$$
z\circ u=b^j\in \cA.
$$
As we did at the beginning of the proof, it follows that $\lg b\rg$ is an $\cA$-subgroup. Thus the basic set containing $b$ could not be as $D$ in this subcase, a contradiction. Similarly, any other singleton support of $x\circ y$ could not be a basic set.  

\medskip

Hence we have
$$
\{a^{i}b^{j},a^{-i}b^{-j}\}\in  \cD(\cA), \, {\rm or} \, \{a^ib^j,a^{nj-i}b^{j}\}\in  \cD(\cA),
$$
 $$
 {\rm or}, \, \{a^ib^j, a^{i-nj}b^{-j}  \}\in \cD(\cA), \, {\rm or}\, \{a^{i}b^{j},a^{nj-i}b^j,a^{-i}b^{-j},\\a^{i-nj}b^{-j}\}\in  \cD(\cA).
 $$

\medskip

{\bf Claim 4.}\, {\it If $2i\ne nj$ and $i,j$ are nonzero, then $\{a^{i}b^{j},a^{-i}b^{-j}\}$ can not be a basic set.}

\proof By way of contradiction, assume that~$\{a^{i}b^{j},a^{-i}b^{-j}\}$~is a basic set. ~Then
 $$
 \{a^{i-nj}b^{-j},a^{nj-i}b^{j}\}
 $$
 is also a basic set.

 Note that
 $$
 v:=a^{2i}b^{j}+b^{j}+b^{-j}+a^{-2i}b^{-j}=(a^{i}+a^{-i})(a^{i}b^{j}+a^{-i}b^{-j})
 \in \cA.
$$
It follows that
$$
u\circ v=b^j+b^{-j}\in \cA.
$$
Consequently, $\lg \{b^j, b^{-j}\}\rg$ is an $\cA$-subgroup, and so is $\lg b \rg$, a contradiction.
\eprf

\medskip

{\bf Claim 5.}\, {\it If $2i\ne nj$ and $i,j$ are nonzero, then $\{a^ib^j,a^{nj-i}b^{j}\}$ can not be a basic set.}

\proof By way also of contradiction, assume that~$\{a^ib^j,a^{nj-i}b^{j}\}$~is a basic set.
Then
$$
w:=(a^i+a^{-i})(a^ib^j+a^{nj-i}b^j)=a^{2i}b^j+a^{nj}b^j+b^j+a^{nj-2i}b^j\in \cA.
$$
Then
$$
w\circ u=b^j+a^{nj}b^j\in \cA.
$$
Thus,
$$
(w\circ u)\circ \underline{D^{(j-1)}}\in \cA,
$$
which yields that
$$
w':= a^{n(2j-1)}b^{2j-1}+a^nb+a^{nj}b^{2j-1}+a^{nj}b+a^{n(j-1)}b^{2j-1}+a^{n(-j+1)}b+b^{2j-1}+b\in \cA.
$$
Since $j\ne 0$, it then follows that
$$
a^nb+b=w'\circ \underline{D}\in \cA,
$$
a contradiction. \eprf
\medskip

{\bf Claim 6.}\, {\it If $2i\ne nj$ and $i,j$ are nonzero, then $\{a^ib^j,a^{i-nj}b^{-j}\}$ can not be a basic set.}

\proof If the claim is false, then
$$
(a^ib^j+a^{i-nj}b^{-j})^2=a^{2i}b^{2j}+a^{2i-2nj}b^{-2j}+2a^{2i-nj}\in \cA,
$$
which implies that $\{a^{2i-nj}\}\in \cD(\cA)$ by Proposition \eqref{913a}. This is a contradiction.
\eprf

\medskip

Consequently, we have that
$$
\{a^{i}b^{j},a^{nj-i}b^j,a^{-i}b^{-j},\\a^{i-nj}b^{-j}\}\in  \cD(\cA), \, \forall i,j\in \Z.
$$
Observe that there exists  $\varphi_2\in \aut(G)$ satisfying
$$
\varphi_2(a)=a^{-1}, \varphi_2(b)=a^nb.
$$
In addition, $G$ has an automorphism $\iota$ such that
$$
\rho(a)=a^{-1}, \rho(b)=b^{-1}.
$$
It is easy to check that $\cA$ is the orbit Schur ring associated with $\lg \varphi_2, \iota\rg$, which is case $(vii)$ in Theorem \ref{1451b}.

\medskip

To complete the proof of Theorem \ref{1451b}, it suffices to show that Schur rings of type $(viii)$ in the theorem is not traditional.

\medskip

Let $\cA$ be a Schur ring over $G$ of type $(viii)$ in Theorem \ref{1451b}. Obviously, $\cA$ can not be the trivial Schur ring. Since $\cZ\times\cZ$ does not have nonidentity finite subgroups, $\cA$ can not be a wedge product Schur ring.

\medskip

Assume that $\cA$ is an orbit Schur ring. Since each basic set in $\cD(\cA)$ has size at most $2$, there exists $\rho\in \aut(G)$ such that
\qtnl{1011a}
\rho(a)=a^{-1},\rho(b)=a^nb^{-1}, \qaq \rho(ab)=a^{n+1}b^{-1},
\eqtn
as $\{a, a^{-1}\}$, $\{b, a^nb^{-1}\}$, and $\{ab, a^{n+1}b^{-1}\}$ are basic sets.

\medskip

However, by the first two equalities in  \eqref{1011a}, one can see that
$$
\rho(ab)=\rho(a)\rho(b)=a^{n-1}b^{-1},
$$
which contadicts the third equality in \eqref{1011a}.

\medskip

Assume that $\cA$ is of tensor product type. Then there exist nontrivial $\cA$-subgroups $H,K$ such that
$$
G=H\times K, \qaq \cA=\cA_H\otimes \cA_K.
$$

Obviously, there exist $D_1\in \cD(\cA_H)$ with $|D_1|=2$ and $D_2\in \cD(\cA_K)$ with $|D_2|=2$. Thus, as a basic set of $\cA$, $D_1\times D_2$ would have size $4$, a contradiction.

\medskip

The proof of Theorem \ref{1451b} is complete.

\eprf

\medskip

\section{Acknowledgment}

The work of the first author is supported by Natural Science Foundation of China (No. 12371019, No. 12161035.)

\end{document}